\documentclass[12pt]{article}

\usepackage{fullpage}

\usepackage{mathrsfs}
\usepackage{comment}
\usepackage{centernot}
\usepackage{mathtools}
\usepackage{stmaryrd}

\usepackage{latexsym,amssymb,amstext,amsmath}
\usepackage{hyperref}
\usepackage{graphbox,graphicx}
\usepackage{subcaption}
\usepackage{amsthm}
\usepackage{esint}
\usepackage{float}
\usepackage{color} 
\usepackage{tikz}

\numberwithin{equation}{section}

\def\be {\begin{equation}}
\def\ee {\end{equation}}
\def\bea {\begin{eqnarray}}
\def\eea {\end{eqnarray}}

\def\M {\mathcal{M}}

\newtheorem{theorem}{Theorem}[section]

\newtheorem{proposition}[theorem]{Proposition}

\theoremstyle{definition}

\newtheorem{remark}[theorem]{Remark}

\begin{document}

\title{\textbf{\LARGE  
Factorisation of symmetric matrices and applications in gravitational theories
}}
\author{M.~Cristina C\^amara and Gabriel Lopes Cardoso}
\date{\small 
\vspace{-5ex}
\begin{quote}
\emph{
\begin{itemize}
\item[]
Center for Mathematical Analysis, Geometry and Dynamical Systems,\\
  Department of Mathematics, 
  Instituto Superior T\'ecnico, Universidade de Lisboa,\\
  Av. Rovisco Pais, 1049-001 Lisboa, Portugal
  \end{itemize}
}
\end{quote}
{\tt 
cristina.camara@tecnico.ulisboa.pt, gabriel.lopes.cardoso@tecnico.ulisboa.pt}
}
\maketitle

\centerline{\it To Yuri Karlovich on his 75th birthday}

\begin{abstract}
\noindent
We consider the canonical Wiener-Hopf factorisation of $2 \times 2$ symmetric matrices $\mathcal M$ with respect to a contour $\Gamma$. For the case that the quotient $q$ of the two diagonal elements of $\mathcal M$ is a rational function, 
we show that due to the symmetric nature of the matrix $\M$, the second column
in each of the two matrix factors that arise in the factorisation is determined in terms of the first column in each of these matrix factors, by multiplication by a rational matrix, and we give a method for determining the second columns of these factors. We illustrate our method with two examples 
in the context of a Riemann-Hilbert approach to obtaining solutions to the Einstein field equations.
\\

%
%
\end{abstract}

\section{Introduction}

Let $\Gamma$ be a simple closed contour in the complex plane encircling the origin and denote by $\mathbb{D}^+_{\Gamma}$ and $\mathbb{D}^-_{\Gamma}$ the interior and the exterior of $\Gamma$.

If $\M$ is an $n \times n$ matrix function whose elements are in $L^{\infty} (\Gamma)$, i.e., essentially bounded functions on $\Gamma$, a bounded Wiener-Hopf (WH for short) factorisation of $\M$ with respect to (w.r.t.) $\Gamma$ is a representation
\bea
\M (\tau) = \M_- (\tau) \, D(\tau) \, \M_+ (\tau) \;\;\;,\;\;\; \tau \in \Gamma \;,
\label{MMDM}
\eea
where $\M_{\pm}$ and their inverses $\M_{\pm}^{-1}$ are analytic and bounded in $\mathbb{D}^{\pm}_{\Gamma}$ (we say that their elements are in $H_{\pm}^{\infty}$, respectively), 
and 
$D(\tau) = \text{diag} ( \tau^{k_1}, \tau^{k_2}, \dots, \tau^{k_n})$ with $k_1 \leq k_2 \leq \dots \leq k_n, \; k_i \in \mathbb{Z} $ for $i=1, 2, \dots, n$. If $k_i = 0$ for all $i=1,2, \dots, n$, then we say that
\bea
\M (\tau) = \M_- (\tau) \, \M_+ (\tau) \;\;\;,\;\;\; \tau \in \Gamma \;,
\eea
is a canonical bounded WH factorisation. In what follows, we will omit the term 'bounded'.

Denoting by $C^{\mu}$  the algebra of all H\"older continuous functions with exponent $\mu \in \, ]0,1[$ 
defined on $\Gamma$ \cite{MP}, 
if $\M$ is invertible in $(C^{\mu})^{n \times n}$, i.e. $\M \in (C^{\mu})^{n \times n}$ and $\det \M \neq 0$ on $\Gamma$, then $\M$ admits a factorisation of the form \eqref{MMDM} with $\M_{\pm} \in 
(C^{\mu}_{\pm})^{n \times n}$, where $C^{\mu}_{\pm} = C^{\mu} \cap H_{\pm}^{\infty}$ \cite{MP}. If the factorisation is canonical, then it is unique if we impose a normalising condition on one of the factors; here we will look for canonical WH factorisations with the factor $\M_+$ normalised to the identity at $0$, in which case it will be denoted by $X$, 
\bea
\M = \M_- \, X  \;\;\; \text{on} \;\;\; \Gamma, \;\;\; \text{with} \; X(0) = \mathbb{I} \;.
\label{WHMf}
\eea

We will be particularly interested in $2 \times 2$ matrix functions $\M$, having in view applications of WH factorisation to solving certain gravitational field equations. For these applications, $\M$ must be a symmetric $2 \times 2$ matrix function
\bea
\M = \begin{bmatrix}
    a & b \\
    b & d
\end{bmatrix} 
\eea
with $a/d \in \mathcal{R}$, where $\mathcal{R}$ denotes the space of all rational functions without poles on $\Gamma$, i.e.
\bea
\frac{a}{d} = \frac{p_1}{p_2} =: q \;,
\label{q}
\eea
where $p_1$ and $p_2$ are polynomials and $p_2$ does not vanish on $\Gamma$, and $\M$ must admit a canonical WH factorisation.

One may then naturally ask how all these conditions, in particular
that $\M$ is symmetric, are reflected in the form of the factors $\M_{\pm}$ of a canonical WH factorisation.

We study this question in Section \ref{sec:meth} by applying and extending an approach which was first presented in \cite{CSC}, although in a different form, and we show that there indeed exists a certain relation between the two columns in each of the factors $\M_-$ and $X^{-1}$ (related by $\M X^{-1} = \M_-$), which is determined
by the structure of the original matrix $\M$ in terms of its symmetry
and the quotient $q$ between its diagonal elements, in such a way that the second column can be obtained from the first by multiplication by a certain rational matrix.
This, on the one hand,  brings out a connection between the form of the matrix $\M$ that is to be factorised and a certain structure of the factors $\M_{\pm}$, which  we present in Theorem \ref{theo:r1r2}; on the other hand, it may also allow for a nicer determination of one of the columns in the factors by avoiding repetition of similar calculations, as we illustrate in the examples 
in Section \ref{sec:EFE}.

The question of how to determine the second columns in a canonical WH factorisation from the first columns has also been studied for other classes of $2 \times 2$ matrix functions, see for instance \cite{CS,BASTOS2003347}.

This also naturally leads to the question of how to determine the first column in $X^{-1}$ and in $\M_-$.
Here we focus on rational matrices possessing a canonical WH factorisation, which are of great importance when considering applications in gravitational theories. Obtaining the first columns is equivalent to determining the (unique) solution to 
\bea
\M \phi_+ = \phi_- \;\;\; \text{on} \;\;\; \Gamma, \;\;\; \text{with} \;\;\; \phi_{\pm} = (\phi_{1 \pm}, \phi_{2 \pm}) \in \left( C^{\mu}_{\pm}  \right)^2 \;\;\;,\;\;
\phi_+ (0)=(1,0) \;.
\eea
We show that it is also possible to simplify the calculations for the first column, reducing the problem of analyticity of the solution to  that of the first component $\phi_{1 \pm}$ (in $\mathbb{D}^{\pm}_{\Gamma}$).

It is well known that Wiener-Hopf factorisation is very important in the study of singular integral equations, convolution equations and in many applications in Mathematics, Physics and Engineering \cite{MP,Y}.
Determining explicit solutions to the Einstein field equations by means of a Riemann-Hilbert (RH) approach is one of the most recent applications of WH factorisation theory.

We illustrate these results in Section \ref{sec:EFE} by applying them to solving
the Einstein field equations by a Riemann-Hilbert approach based on 
\cite{Breitenlohner:1986um,Aniceto:2019rhg}.


\section{Canonical factorisation and structure of the factors \label{sec:meth}}

In what follows, for any algebra $A$, let ${\cal G} \, A$ denote the group of invertible elements in $A$.

Let 
\bea
\M (\tau) = \begin{bmatrix}
    a(\tau) & b(\tau) \\
    b(\tau) & d(\tau)
\end{bmatrix}
\eea
with $a, b, d $ analytic in a neighbourhood of $\Gamma$, where $\Gamma$ is a simple closed contour, be such that $\M$ admits a canonical WH factorisation 
w.r.t. $\Gamma$
of the form \eqref{WHMf}, 
\bea
\M = \M_- \, X \;\;\; \text{on} \;\;\; \Gamma, \;\;\; \text{with} \; X(0) = \mathbb{I} \;.
\eea
Then, $d = \det \M$ also admits a canonical WH factorisation w.r.t. $\Gamma$, 
\bea
d = d_- \, d_+  \;\;\; \text{with} \;\;\; d_{\pm} \in \mathcal{G} C^{\mu}_{\pm}  \;\;\;,\;\;\; d_+(0) = 1 \;.
\label{dd+d-}
\eea
Let 
\bea
\frac{a}{d}= q \in  \mathcal{R} \;\;\;,\;\;\; q = \frac{p_1}{p_2} \;,
\label{adq}
\eea
where $p_1$ and $p_2$ are polynomials without common zeroes 
and such that $q$ is bounded at $\infty$.

The symmetric structure of $\M$, which will be reflected in the form of the factors $X^{-1}$ and $\M_-$, can be characterised by the following relation.

{\proposition 
Let $q$ be defined by \eqref{adq} and let
\bea
Q_1 = {\rm diag} \, ( 1 \; -q) \;\;\;,\;\;\; Q_2 = {\rm diag} \, ( q \; -1) \;.
\label{Q1Q2diag}
\eea
Then
\bea
\M Q_1 \M = d \, Q_2 \;\;\; \text{with} \;\;\; d = \det \M = a d - b^2 \;.
\label{MQ1MQ2}
\eea

}

{\remark
Note that, conversely, if $\tilde \M$ is a $2 \times 2$ matrix
\bea
{\tilde \M} = \begin{bmatrix}
    a & b\\
    c & d 
    \end{bmatrix} \;\;\; \text{with} \;\;\; a, d, ad-bc \neq 0 \, 
\eea
and
\bea
{\tilde \M}^T Q_1 {\tilde \M} = (\det {\tilde \M}) \, Q_2
\eea
with $Q_1 = {\rm diag} \, ( 1 \; -q) \;,\; Q_2 = {\rm diag} \, ( q \; -1)$ for some $q$, then we have that ${\tilde \M}^T = {\tilde \M}$ and $q = a/d$. So the relation \eqref{MQ1MQ2} does indeed characterise the structure of $\M$.\\

}

The first columns of $X^{-1}$ and $\M_-$, denoted $f_+$ and $f_-$ respectively, are
the unique solution to 
\bea
\M f_+ = f_- \,\;\; \text{with} \;\;\; f_{\pm} \in \left( C^{\mu}_{\pm}  \right)^2 \;\;\;,\;\;\; f_+ (0) = 
\begin{bmatrix}
    1\\ 0
\end{bmatrix} \;;
\label{Mf+f-}
\eea
the second columns of $X^{-1}$ and $\M_-$, denoted $s_+$ and $s_-$, respectively, constitute the unique solution
\bea
\M s_+ = s_- \,\;\; \text{with} \;\;\; s_{\pm} \in \left( C^{\mu}_{\pm}  \right)^2 \;\;\;,\;\;\; s_+ (0) = 
\begin{bmatrix}
    0\\ 1
\end{bmatrix} \;.
\label{Ms+s-}
\eea
We therefore have that 
\bea
\M \, \begin{bmatrix}
f_+ & s_+
\end{bmatrix} = \begin{bmatrix}
f_- & s_-
\end{bmatrix} \;,
\label{Mfsfs}
\eea
where
\bea
X^{-1} =  \, \begin{bmatrix}
f_+ & s_+
\end{bmatrix} \;\;\;,\;\;\; \M_- =  \, \begin{bmatrix}
f_- & s_-
\end{bmatrix} \;\;\;,\;\;\; X(0) = \mathbb{I} \;.
\eea

To study the relation between $f_{\pm}$ and $s_{\pm}$ in this case, we will use the following results, which can be easily verified.

{\proposition \label{prop:JAJ}
If $A$ is a $2 \times 2$ matrix, $A = [ f \; s]$, and
\bea
J = \begin{bmatrix}
0 & -1 \\
1 & 0 
\end{bmatrix} \;,
\eea
then 

(i) $A J A^T = (\det A) \, J$, 

(ii) $\det A = s^T J f $, 

(iii) $A^{-1} = - ( \det A)^{-1} \, J A^T J $ if $\det A \neq 0$. \\}

Also note that, for any vector $f$ with two components,
\bea
f^T J f = 0 \;.
\label{fjf}
\eea

{\proposition \label{prop:QQJJff}
Let $\M f_+ = f_-$ with $\M$ satisfying \eqref{MQ1MQ2}. Then 
\bea
\M \left( J Q_2 f_+ \right) = J Q_1 f_- \;.
\label{MJMJQQ}
\eea

}

\begin{proof}

We have that $Q_2 = d^{-1} \M Q_1 \M$, hence
\bea
\M \left( J Q_2 f_+ \right) = d^{-1} \M J \left( \M Q_1 \M \right) f_+ = \left(
d^{-1} \M J \M \right) Q_1 \left( \M f_+ \right) = J Q_1 f_- \;,
\eea
where we used Proposition \ref{prop:JAJ} (i).

\end{proof}

Thus we have
\bea
\M \begin{bmatrix}
f_+ & J Q_2 f_+ 
\end{bmatrix} = \begin{bmatrix} f_- & J Q_1 f_-
\end{bmatrix}
\eea
and, applying determinants and using \eqref{dd+d-}, we obtain
\bea
d_+ \det \begin{bmatrix}
f_+ & J Q_2 f_+ 
\end{bmatrix} = d^{-1}_- \det \begin{bmatrix} f_- & J Q_1 f_-
\end{bmatrix} = r_1 \;,
\label{r1dd}
\eea
where $r_1 \in \mathcal{R}$, with $r_1$ bounded at $\infty$, since
the left hand side of the first equality is meromorphic in 
$\mathbb{D}^+_{\Gamma}$ and the right hand side is meromorphic in 
$\mathbb{D}^-_{\Gamma}$ and bounded at $\infty$. Moreover, from 
\eqref{Q1Q2diag} we see that the poles of $r_1$ must be those of $q$ (see \eqref{r1dd}), so we may assume that
\bea
r_1 = \frac{{\tilde p}_1}{p_2} \;,
\label{r1tp1p2}
\eea
where ${\tilde p}_1$ is a polynomial of degree not greater than
$\deg (p_2)$. We may also assume that $r_1$ does not have zeroes on $\Gamma$, since, given the analyticity of $a, b, d$ in a neighbourhood of $\Gamma$, the latter can be deformed if necessary.

Now, from \eqref{r1dd}, we have 
\bea
d_+ \det \begin{bmatrix}
f_+ & r_1^{-1} J Q_2 f_+ 
\end{bmatrix} = d^{-1}_- \det \begin{bmatrix} f_- & r_1^{-1} J Q_1 f_-
\end{bmatrix} = 1 \;.
\label{r1dd1}
\eea
On the other hand, 
applying determinants to \eqref{Mfsfs}, we obtain
\bea
d_+ \, \det \begin{bmatrix}
f_+ & s_+
\end{bmatrix} = d_-^{-1} \,  \det \begin{bmatrix}
f_- & s_-
\end{bmatrix} 
\label{dfsdfs}
\eea
and, 
since the left hand side of 
\eqref{dfsdfs} 
is in $C^{\mu}_+ $, while the right hand side is in 
$C^{\mu}_- $, we conclude that both sides are equal to a constant $k=1$, taking into account that the left hand side must equal $1$ at $\tau =0$. Thus we have
\bea
d_+ \, \det \begin{bmatrix}
f_+ & s_+
\end{bmatrix} = 1 = d_-^{-1} \, \det \begin{bmatrix}
f_- & s_-
\end{bmatrix} \;.
\label{dfsdfs2}
\eea
Then, combining \eqref{dfsdfs2} with \eqref{r1dd1} gives
\bea
d_+ \det \begin{bmatrix}
f_+ & \;\; s_+ - r_1^{-1} J Q_2 f_+ 
\end{bmatrix} = d^{-1}_- \det \begin{bmatrix} f_- & \;\; s_- - r_1^{-1} J Q_1 f_-
\end{bmatrix} = 0\;.
\label{r1dd2}
\eea
It follows that
\bea
s_+ - r_1^{-1} J Q_2 f_+ &=& \lambda_1 \, f_+ \;, 
\label{spfpQQ}\\
s_- - r_1^{-1} J Q_1 f_- &=& \lambda_2 \, f_- \;,
\label{smfmQQ}
\eea
where $\lambda_1$ and $\lambda_2$ are functions of $\tau$ which, due to \eqref{Mf+f-},  \eqref{Ms+s-} and \eqref{MJMJQQ}
must satisfy
\bea
\lambda_1 = \lambda_2 =: \lambda \;.
\eea

{\proposition

With the notation above, $\lambda \in \mathcal{R}$ and its poles are the zeroes of $r_1$.
}

\begin{proof}

{From} \eqref{fjf} and \eqref{spfpQQ}, with $\lambda_1 = \lambda$,
we have that
\bea
0 &=& s_+^T J s_+ = \lambda s_+^T J f_+ + r_1^{-1} s_+^T J J Q_2 f_+
= \lambda \det 
\begin{pmatrix} f_+ & s_+ \end{pmatrix}
-  r_1^{-1} s_+^T Q_2 f_+ \nonumber\\
&=&\lambda \, d^{-1}_+ -  r_1^{-1} s_+^T Q_2 f_+ \;,
\eea
and hence
\bea
\lambda = d_+ \, r_1^{-1} s_+^T Q_2 f_+ \in H_+^{\infty} + \mathcal{R} \;.
\label{r1R1}
\eea
On the other hand, from \eqref{fjf} and \eqref{smfmQQ}
with $\lambda_2 = \lambda$,
we get in an analogous way
\bea
\lambda = d_-^{-1} \, r_1^{-1} s_-^T Q_1 f_- \in H_-^{\infty} + \mathcal{R} \;.
\label{r1R2}
\eea
It follows from \eqref{r1R1} and \eqref{r1R2} that 
$\lambda \in \mathcal{R}$. Since $r_1$ and $q$ have the same denominator, we see that the poles of $Q_1$ and $Q_2$ are cancelled and therefore the poles of $\lambda$ are those of $ r_1^{-1}$, i.e., the zeroes of $r_1$.
\end{proof}

As a consequence of the previous results, we obtain the following relation between the two columns of $X^{-1}$ and $\M_-$.

{\theorem \label{theo:r1r2}
With the same notation as in Proposition \ref{prop:QQJJff}, there exists $r_2 \in \mathcal{R}$, bounded at $\infty$, with the same
poles as $q$, such that
\bea
s_+ &=&  r_1^{-1} \left( r_2 \mathbb{I} + J Q_2 \right)  f_+ 
\Longleftrightarrow
\left\{ 
\begin{matrix}
s_{1+} &=&  r_1^{-1} \left( r_2 f_{1+} + f_{2+}   \right) \\
s_{2+} &=& r_1^{-1} \left( r_2 f_{2+} + q f_{1+}   \right) 
\end{matrix}
\right.
\label{sr2pfpQQ}\\
s_- &=& r_1^{-1} \left(r_2 \mathbb{I} +  J Q_1 \right)  f_- 
\Longleftrightarrow
\left\{ 
\begin{matrix}
s_{1-} &=&  r_1^{-1} \left( r_2 f_{1-} + q f_{2-}   \right) \\
s_{2-} &=& r_1^{-1} \left(  r_2 f_{2-} + f_{1-}   \right) 
\end{matrix}
\right.
\label{sr1mfmQQ}
\eea

}

\begin{proof}
    
{From} \eqref{spfpQQ} and \eqref{smfmQQ} we have 
\bea
s_+ &=&  r_1^{-1} \left( \lambda r_1 \mathbb{I} + J Q_2 \right)  f_+ \Longleftrightarrow
\left\{ 
\begin{matrix}
s_{1+} &=&  r_1^{-1} \left( \lambda r_1 f_{1+} + f_{2+}   \right) \\
s_{2+} &=& r_1^{-1} \left( \lambda r_1 f_{2+} + q f_{1+}   \right) 
\end{matrix}
\right.
\label{sr2pfpQQ2}\\
s_- &=& r_1^{-1} \left(\lambda r_1 \mathbb{I} +  J Q_1 \right)  f_- \Longleftrightarrow
\left\{ 
\begin{matrix}
s_{1-} &=&  r_1^{-1} \left( \lambda r_1 f_{1-} + q f_{2-}   \right) \\
s_{2-} &=& r_1^{-1} \left( \lambda r_1 f_{2-} + f_{1-}   \right) 
\end{matrix}
\right.
\label{sr1mfmQQ2}
\eea
where
\bea
\lambda r_1 = d_-^{-1} s_-^T Q_1 f_- = d_+ s_+^T Q_2 f_+ \;,
\eea
so the result holds with $r_2 = \lambda r_1$. Note that $r_2$ is bounded at $\infty$, since $q$ is bounded there.

\end{proof}

Now we address the question of determining $r_2$. It must be such that $s_{\pm}$, given by \eqref{sr2pfpQQ} and \eqref{sr1mfmQQ}, are analytic in $\mathbb{D}^{\pm}_{\Gamma} $, 
$s_{\pm} \in (H_{\pm}^{\infty})^2$,
and $s_+(0) = \begin{bmatrix} 0 \\ 1 \end{bmatrix}$; in that case
$s_+$ and $s_-$ are the second columns of $X^{-1}$ and $\M_-$, respectively, as formulated next.

{\proposition

If $r_2 \in \mathcal{R}$ is such that $s_+$ and $s_-$, defined by the right hand side 
of \eqref{sr2pfpQQ} and \eqref{sr1mfmQQ}, respectively, are functions in 
$(H_{\pm}^{\infty})^2$, with $s_+(0) = \begin{bmatrix} 0 \\ 1 \end{bmatrix}$, then $s_+$ and $s_-$ are the second columns of 
$X^{-1}$ and $\M_-$, respectively. In that case, $r_2$
is unique.
}

\begin{proof}

We have that
\bea
\M s_+ = \M \left(  r_1^{-1} \left( r_2 \mathbb{I} + J Q_2 \right)  f_+ \right) =  r_1^{-1} \left( \M r_2 f_+ + \M J Q_2  f_+ \right) = r_1^{-1} \left(  r_2 f_- + J Q_1  f_- \right) = s_-
\;, \nonumber
\eea
and since  $s_{\pm} \in (H_{\pm}^{\infty})^2$ with $s_+(0) = \begin{bmatrix} 0 \\ 1 \end{bmatrix}$, we have the second columns of $X^{-1}$ and $\M_-$.
Moreover, from \eqref{sr2pfpQQ} we get
\bea
r_1 s_+^T J s_+ = r_2 s_+^T J f_+ - s_+^T Q_2  f_+
\Longleftrightarrow
0 = r_2 \, d_+^{-1} - s_+^T Q_2  f_+ \Longleftrightarrow r_2 = d_+
\, s_+^T Q_2  f_+ \;, \nonumber
\eea
and the uniqueness of $r_2$ follows from the uniqueness of $f_+$ and $s_+$.

\end{proof}

{\remark
Note that $r_2$ can also be obtained from \eqref{sr1mfmQQ} and, in that case, one obtains $r_2 = d_-^{-1}
\, s_-^T Q_1  f_- $. Indeed, from $\M s_+ = s_-$ we get that
$s_+^T \M = s_-^T$ and, from \eqref{MQ1MQ2},
\bea
s_+^T \M Q_1 \M f_+ = s_-^T Q_1 f_- \Longleftrightarrow 
s_+^T d \, Q_2 f_+ = s_-^T Q_1 f_- \Longleftrightarrow 
d_+ \, s_+^T Q_2 f_+ = d_-^{-1} \, s_-^T Q_1 f_- \;. \nonumber
\eea
\\

}

Since $r_2$ has the form 
\bea
r_2 = \frac{R_2}{p_2} \;,
\eea
where $p_2$ is the denominator of $q$ as in \eqref{q} and $R_2$ is a polynomial of degree not greater
than $\deg (p_2)$ (since $r_2$ is bounded at $\infty$), we only have to determine the coefficients of $R_2$, which must be such that the zeroes
of $(r_2 \mathbb{I} + J Q_2) f_+ $ cancel the zeroes of $r_1$ in $\mathbb{D}^+_{\Gamma}$ and the zeroes of
$(r_2 \mathbb{I} + J Q_1) f_- $ cancel the zeroes of $r_1$ in $\mathbb{D}^-_{\Gamma}$. 
We now present a systematic method to obtain those coefficients in the case where $r_1$ has simple or double zeroes; the method can however be generalised, following the same reasoning, for higher order of zeroes.

Let each zero of $r_1$ be denoted by $z_i^+$ if it belongs to $\mathbb{D}^+_{\Gamma}$, and by 
 $z_i^-$ if it belongs to $\mathbb{D}^-_{\Gamma}$. Imposing that $(r_2 \mathbb{I} + J Q_2) f_+ $ vanishes
 at a zero $z_i^+$ of $r_1$ means that, for $f_+ = (f_{1+},  f_{2+})$, and $r_1$ given by \eqref{r1tp1p2}, we must have
 \bea
 \left\{
\begin{matrix}
 R_2 (z_i^+) \, f_{1 +} (z_i^+) + p_2(z_i^+) \, f_{2 +}  (z_i^+) &=& 0 \;, 
 \label{rzf1} \\
R_2 (z_i^+) \,  f_{2 +}  (z_i^+) + p_1(z_i^+)  \, f_{1+} (z_i^+)  &=& 0 \;.
\end{matrix}
\right.
\label{r2zf2}
 \eea
Note that, since $r_1$ is given by \eqref{r1dd}, 
we have that 
\bea
r_1 = d_+ \, \frac{p_1 \, f_{1+}^2 - p_2 \, f_{2+}^2}{p_2} \;,
\eea
so 
\bea
p_1 (z_i^+) \, f_{1+}^2 (z_i^+) = p_2 (z_i^+)\, f_{2+}^2 (z_i^+) \;,
\label{p1p2ffzp}
\eea
and it follows that, since $ f_{1+}$ and  $f_{2+}$ cannot vanish simultaneously,
\bea
f_{2+} (z_i^+) = 0 \Longrightarrow p_1 (z_i^+) = 0 \;, \\
f_{1+} (z_i^+) = 0 \Longrightarrow p_2 (z_i^+) = 0 \;.
\eea
We will show that \eqref{r2zf2} reduces to just one equation, for which we consider three cases:

(i) if $f_{2+} (z_i^+) = 0 $, the first equation in \eqref{r2zf2} is equivalent to
\bea
R_2 (z_i^+) = 0 \;,
\label{R2z0}
\eea
while the second equation is satisfied for any $R_2 (z_i^+)$;

(ii) if $f_{1+} (z_i^+) = 0 $, the second equation in \eqref{r2zf2} is equivalent to \eqref{R2z0},
while the first equation is satisfied for any $R_2 (z_i^+)$;

(iii) if $f_{1+} (z_i^+) , f_{2+} (z_i^+) \neq  0 $, then multiplying the second equation in \eqref{r2zf2} 
by $f_{1+} (z_i^+)$ we get
\bea
R_2 (z_i^+) \,  f_{2 +}  (z_i^+) \, f_{1 +}  (z_i^+) + p_1(z_i^+)  \, f^2_{1+} (z_i^+)  &=& 0 \nonumber\\
\Longleftrightarrow 
R_2 (z_i^+) \,  f_{2 +}  (z_i^+) \, f_{1 +}  (z_i^+) + p_2 (z_i^+)\, f_{2+}^2 (z_i^+)  &=& 0 \nonumber\\
\Longleftrightarrow 
f_{2 +}  (z_i^+) \left[R_2 (z_i^+)  \, f_{1 +}  (z_i^+) + p_2 (z_i^+)\, f_{2+} (z_i^+) \right]  &=& 0 \;,
\eea
which is equivalent to the first equation in  \eqref{r2zf2}, i.e.
\bea
R_2 (z_i^+)  \, f_{1 +}  (z_i^+) + p_2 (z_i^+)\, f_{2+} (z_i^+) =0 \;.
\label{R2ffz0}
\eea

Analogously, for $f_- = (f_{1-},  f_{2-})$ and a zero $z_i^-$  of $r_1$, we get
\bea
R_2 (z_i^-)= 0 \;\;\;,\;\;\; \text{if} \;\;\;\; f_{1-} (z_i^-) = 0  \;\;\; \text{or} \;\;\;  f_{2-} (z_i^-) = 0 \;,
\label{R2zm0}
\\
R_2 (z_i^-)  \, f_{1 -}  (z_i^-) + p_1 (z_i^-)\, f_{2-} (z_i^-) =0 \;\;\;,\;\;\; 
 \text{if} \;\;\;\; f_{1-} (z_i^-), f_{2-} (z_i^-)  \neq 0 \;.
 \label{R2ffzm0}
\eea

If all the zeroes of $r_1$ are simple zeroes, then \eqref{R2z0} / \eqref{R2ffz0} and \eqref{R2zm0} / \eqref{R2ffzm0}
provide a system of equations allowing to determine all but one coefficient of $R_2$; the remaining coefficient is obtained from the normalising condition $s_{1+} (0) = 0$.

{\remark \label{normcond}

Note that the normalising condition $s_{1+} (0) = 0 $ implies that $s_{2+} (0) = 1 $, as follows.
First we note that
\bea
r_1 = 
d_+ \, \left( q \, f_{1 +}^2 - f_{2 +}^2 \right) \;.
\label{r1qf1f2}
\eea
Using \eqref{sr2pfpQQ2} and \eqref{r1qf1f2} we have
\bea
s_{2+} f_{1+} &=& r_1^{-1} \left( r_2 \, f_{1+} f_{2+} + q\, f^2_{1+} \right) \nonumber\\
&=&  r_1^{-1} r_2 \, f_{1+} f_{2+} + r_1^{-1} \left( d_+^{-1}  r_1 + f^2_{2+} \right) \nonumber\\
&=& r_1^{-1} \left( r_2 \, f_{1+} + f_{2+} \right) f_{2+} +  d_+^{-1} \nonumber\\
&=& s_{1+} f_{2+} + d_+ \;.
\eea
Imposing $s_{1+} (0) = 0 $ and using $d_+(0) = 1$, this yields $s_{2+} (0) = 1 $.
\\

To extend the method presented above to the case where $r_1$ has double zeroes, let us assume that 
${\tilde p}_1$ has a double zero at 
$z_i^+ \in \mathbb{D}^+_{\Gamma}$. 
Using the above results, we see that the conditions that we have to impose on $r_2$ are the same as above and moreover
 \bea
 \left\{
 \begin{matrix}
   R'_2 (z_i^+) \, f_{1 +} (z_i^+) +  R_2 (z_i^+) \, f'_{1 +} (z_i^+) +  p'_2 (z_i^+) \, f_{2 +} (z_i^+)
+ p_2 (z_i^+) \,   f'_{2 +}  (z_i^+) &=& 0 \;, \\
  R'_2 (z_i^+) \,  f_{2 +}  (z_i^+) + 
  R_2 (z_i^+) \,  f'_{2 +}  (z_i^+) +
   p_1'(z_i^+)  \, f_{1+} (z_i^+) + p_1(z_i^+)  \, f'_{1+} (z_i^+) &=& 0 \;.
   \end{matrix}
   \right.
   \label{condr2}
\eea
To show that this pair of conditions can be reduced to one equivalent condition, we consider once again three cases. First, however, note that saying that ${\tilde p}_1$ has a double zero at $z_i^+$ means that
(cf. \eqref{p1p2ffzp})
\bea
p_1'(z_i^+) \, f^2_{1+} (z_i^+) + 2 p_1 (z_i^+) \, f'_{1+} (z_i^+) \,f_{1+} (z_i^+) - 
p_2'(z_i^+) \, f^2_{2+} (z_i^+)  -  2 p_2 (z_i^+) \,f'_{2+} (z_i^+) \,f_{2+} (z_i^+) = 0 \;,
\nonumber
\eea
i.e.
\bea
p_1'(z_i^+) \, f^2_{1+} (z_i^+) +  p_1 (z_i^+) \, f'_{1+} (z_i^+) \,f_{1+} (z_i^+) &=&  
p_2'(z_i^+) \, f^2_{2+} (z_i^+) -  p_1 (z_i^+) \, f'_{1+} (z_i^+) \,f_{1+} (z_i^+) \nonumber\\
&& +  2 p_2 (z_i^+) \,f'_{2+} (z_i^+) \,f_{2+} (z_i^+)  \;.
\label{pffpfp}
\eea
With this in mind,

(i) if $f_{2+} (z_i^+) = 0$, from \eqref{R2z0} and the first equation in \eqref{condr2} we get
\bea
\left\{
\begin{matrix}
    && R_2 (z_i^+) = 0 \\
     && R'_2 (z_i^+) \, f_{1 +} (z_i^+)  
+ p_2 (z_i^+) \,   f'_{2 +}  (z_i^+) =0 \;, \\
\end{matrix}
\right.
\label{f2p0}
\eea
while the second equation in \eqref{condr2} is satisfied for any $R_2$ because 
$f_{2+} (z_i^+) = p_1 (z_i^+) = R_2 (z_i^+) = 0$ and \eqref{pffpfp} implies that $p_1'(z_i^+) =0$;

(ii)  if $f_{1+} (z_i^+) = 0$, the equations analogously reduce to
\bea
\left\{
\begin{matrix}
    && R_2 (z_i^+) = 0 \\
     && R'_2 (z_i^+) \, f_{2 +} (z_i^+)  
+ p_1 (z_i^+) \,   f'_{1 +}  (z_i^+) =0 \;; \\
\end{matrix}
\right.
\eea

(iii) if $ f_{1 +} (z_i^+),  f_{2 +} (z_i^+) \neq 0$, then multiplying the second equation in 
\eqref{condr2} by $f_{1 +}  (z_i^+)$ we see that it is equivalent to 
\bea
  R'_2 (z_i^+) \,  f_{2 +}  (z_i^+) \,  f_{1+} (z_i^+) + 
  p_2 (z_i^+) \,  f_{2 +}  (z_i^+) \,  f'_{2+} (z_i^+) +  p'_2 (z_i^+) \,  f^2_{2 +}  (z_i^+) 
  \nonumber\\
  - 
   p_1(z_i^+)  \, f_{1+} (z_i^+) \,  f'_{1+} (z_i^+) = 0 \;,
   \label{R2pR2}
   \eea
where we used \eqref{R2ffz0} and \eqref{pffpfp}. Now, 
from \eqref{R2ffz0} and \eqref{p1p2ffzp} we have that
\bea
R_2 (z_i^+)  \, f_{2 +}  (z_i^+) + p_1 (z_i^+)\, f_{1+} (z_i^+) =0 \;,
\eea
and substituting in \eqref{R2pR2} we obtain
\bea
 f_{2 +}  (z_i^+) \left[ R'_2 (z_i^+)  \,  f_{1+} (z_i^+) + 
  p_2 (z_i^+) \,   f'_{2+} (z_i^+) +  p'_2 (z_i^+) \,  f_{2 +}  (z_i^+) +  R_2(z_i^+)   \,  f'_{1+} (z_i^+)
  \right] =0 \;, \nonumber
   \eea
which is equivalent to the first condition in \eqref{condr2}. Therefore we find that imposing a double zero for ${\tilde p}_1$ at $z_i^+$ is equivalent to imposing, in this case,
\bea
\left\{
\begin{matrix}
R_2 (z_i^+)  \, f_{1 +}  (z_i^+) + p_2 (z_i^+)\, f_{2+} (z_i^+) =0 \;, \\
  R'_2 (z_i^+) \, f_{1 +} (z_i^+) +  R_2 (z_i^+) \, f'_{1 +} (z_i^+) + \left(  p_2 \, f_{2+} \right)'(z_i^+)  = 0  \;.
\end{matrix}
\right.
\eea

Analogously, for $f_-=(f_{1-}, f_{2-})$ and a double zero $z_i^-$ of $r_1$ in 
$\mathbb{D}^-_{\Gamma}$, we obtain the conditions:

(i) if $f_{1-} (z_i^-) = 0$, 
\bea
\left\{
\begin{matrix}
    && R_2 (z_i^-) = 0 \\
     && R'_2 (z_i^-) \, f_{2-} (z_i^-)  
+ p_2 (z_i^-) \,   f'_{1-}  (z_i^-) =0 \;; \\
\end{matrix}
\right.
\eea

(ii)  if $f_{2-} (z_i^-) = 0$, 
\bea
\left\{
\begin{matrix}
    && R_2 (z_i^-) = 0 \\
     && R'_2 (z_i^-) \, f_{1-} (z_i^-)  
+ p_1 (z_i^-) \,   f'_{2 -}  (z_i^-) =0 \;; \\
\end{matrix}
\right.
\eea

(iii) if $ f_{1 -} (z_i^-),  f_{2 -} (z_i^-) \neq 0$, 
\bea
\left\{
\begin{matrix}
R_2 (z_i^-)  \, f_{1-}  (z_i^-) + p_1 (z_i^-)\, f_{2-} (z_i^-) =0 \;, \\
  R'_2 (z_i^+) \, f_{1-} (z_i^-) +  R_2 (z_i^-) \, f'_{1-} (z_i^-) + \left(  p_1 \, f_{2-} \right)'
  (z_i^-)  = 0  \;.
\end{matrix}
\right.
\label{ffm}
\eea

Proceeding analogously for every double zero of $r_1$, and as in the previous step for every single zero of $r_1$, 
and adding the normalising condition $s_{1+} (0) = 0 $,
we obtain a linear system providing the coefficients of the numerator of $r_2$.

We end this section with a result which can be used to simplify the determination of $f_+$ and $f_-$ above, but can also be extended to more general matrices, when $\M$ is rational. It relates
the two components $\phi_{1+}$ and $\phi_{2+}$ of any solution to the Riemann-Hilbert problem of the form 
\bea
\M \phi_+ = \phi_- \;\;\; \text{on}\;\;\;  \Gamma, \;\;\; \text{with} \;\;\; \phi_{\pm} = \begin{bmatrix} \phi_{1+} \\ \phi_{2+} \end{bmatrix}
\in ( H^{\infty}_{\pm})^2 \;,
\label{RHpphi}
\eea
such as \eqref{Mf+f-} and \eqref{Ms+s-}, in the case that $\M$ is a rational $2 \times 2$ 
(not necessarily symmetric) matrix function
\bea
\M  = \begin{bmatrix}
    q_{11} & q_{12} \\
    q_{21} & q_{22}
\end{bmatrix} \;.
\eea
{From} \eqref{RHpphi} we have that
\bea
q_{11} \phi_{1+} + q_{12} \phi_{2+} = \phi_{1-} = T_1 \in \mathcal{R} \;, \nonumber\\
q_{21} \phi_{1+} + q_{22} \phi_{2+} = \phi_{2-} = T_2 \in \mathcal{R} \;, 
\label{Tsys}
\eea
where $T_1$ and $T_2$ are rational functions bounded at $\infty$, whose denominators are defined by the poles of $q_{ij}$ in $\mathbb{D}^+_{\Gamma}$  ($i,j=1,2$). We can therefore reduce the problem to the following case with polynomial coefficients, 
\bea
p_{11} \phi_{1+} + p_{12} \phi_{2+} = P_1  \;, \nonumber\\
p_{21} \phi_{1+} + p_{22} \phi_{2+} = P_2  \;, 
\eea
where $p_{ij}, P_i$  ($i,j=1,2$) are polynomials. We have the following, which is a slight generalisation of Lemma 3.9
in \cite{Camara:2024ham}.

{\proposition

Let $p_{ij}, P_i$  ($i,j=1,2$) be polynomials such that $p_{11} \, p_{22} - p_{12} \, p_{21}$ does not vanish on $\Gamma$.
Assume moreover that $p_{1i}, p_{2i}$ do not have common zeroes in $\mathbb{D}^+_{\Gamma}$ and consider
the solution 
$\phi = \begin{bmatrix} \phi_1 \\ \phi_2 \end{bmatrix}$ of the system 
 \bea
 \label{ppass}
     \left.
    \begin{array}{lr}
      {p}_{11}  \, \phi_{1} +  {p}_{12} \, \phi_{2} &= P_1 \\
      {p}_{21} \, \phi_{1} +  {p}_{22} \, \phi_{2}  &= P_2
                  \end{array}\right\} \quad \text{on} \quad \Gamma .
                  \label{sysp12}
  \eea
Then, if $\phi_1$ is analytic in $\mathbb{D}^+_{\Gamma}$, $\phi_2$ is also analytic in $\mathbb{D}^+_{\Gamma}$, and vice versa.

}

\begin{proof}

By Cramer's rule we have that
\bea
\phi_1 = \frac{ P_1 \, p_{22} - P_2 \, p_{12}}{p_{11} \, p_{22} - p_{12} \, p_{21}} \;.
\eea
On the other hand, from \eqref{sysp12} we have
\bea
\phi_2 = \frac{P_1 - p_{11} \, \phi_1}{p_{12} } = \frac{ P_2 - p_{21} \, \phi_1}{p_{22}} \;.
\label{phi2p12p22}
\eea
If $\phi_1$ is analytic in $\mathbb{D}^+_{\Gamma}$ and $p_{12}, p_{22}$ do not have common zeroes in $\mathbb{D}^+_{\Gamma}$, then in the neighbourhood of any zero of $p_{12}$ in $\mathbb{D}^+_{\Gamma}$ we see from the second equality in \eqref{phi2p12p22} that $\frac{ P_2 - p_{21} \, \phi_1}{p_{22}}$ must be analytic. Thus $\phi_2$
is analytic in $\mathbb{D}^+_{\Gamma}$.

Conversely, if $\phi_2$ is analytic in $\mathbb{D}^+_{\Gamma}$,
it follows by an analogous argument that $\phi_1$ is also   analytic in $\mathbb{D}^+_{\Gamma}$.

\end{proof}

{\remark

To obtain a solution to \eqref{Mf+f-} when $\M$ is a rational matrix, one has to solve \eqref{Tsys} for the unknowns
$T_1$ and $T_2$, which must be determined such that both $\phi_{1+}$ and $\phi_{2+}$ belong to $H_+^{\infty}$. 
The result given above shows that it suffices to determine $T_1$ and $T_2$ such that $\phi_{1+}$ is 
analytic in $\mathbb{D}^+_{\Gamma}$.

}


\section{The Einstein field equations and the monodromy matrix \label{sec:EFE}}

The Einstein field equations, a system of 10 
nonlinear second order PDE's in 4 variables for the space-time metric $g$,
relate the geometry of space-time to the distribution of matter and energy, described by the 
stress-energy-momentum tensor $T$, within it.
In the following, we will assume the absence of a cosmological constant term in the field equations. 
When $T=0$ in the region under consideration, the field equations are also referred to as the vacuum field equations.

Obtaining exact solutions to the field equations is, in general, a 
non-trivial task. Exact solutions, which play an important role in Physics and Mathematics, can however be obtained
under simplifying assumptions, such as symmetry conditions.
We focus on the subspace of solutions of the field equations possessing two commuting isometries, so that the theory can be reduced to two dimensions using a well-known 2-step procedure (cf. \cite{Lu:2007jc}), and the problem of solving the field equations is reduced to a system of nonlinear second order PDE's depending on two coordinates, which we denote by $\rho$ and $v$, called Weyl coordinates, with $\rho >0$.
We identify these solutions 
with matrix functions $M(\rho, v)$ of class $C^2$, which satisfy the field equations
\bea
d \left( \rho \star A \right) = 0 \;\;\;,\;\;\; A = M^{-1} d M \;,
\label{eomM}
\eea
where $\star$ is the Hodge star operator, $\det M = 1$ and 
$M = M^{\natural}$ (cf. \cite{Lu:2007jc}). Here $\natural$ denotes a certain involution called generalised transposition. When $T=0$, i.e. when dealing with the vacuum field equations, $M$ is a $2 \times 2$ matrix and $\natural$ denotes matrix transposition. 

Determining explicit solutions to the field equations \eqref{eomM} by means of a Riemann-Hilbert (RH) approach is one of the most recent applications of WH factorisation theory.
Here the factorisation is considered with respect to an admissible contour, by which we mean
a closed simple contour in the complex plane, encircling the origin and invariant under the involution $\tau \xmapsto[\, i_{\lambda}]{}
- \frac{\lambda}{\tau}$, where $\lambda = \pm 1$ depending on the physics context. This RH approach
crucially involves introducing a complex parameter $\tau$, called the spectral parameter, which is allowed to vary on an algebraic curve, called the spectral curve, given by the relation 
\bea
\omega = v + \frac{\lambda}{2}    \, \rho \, \frac{\lambda -   \tau^2}{\tau} \;\;\;,\;\;\; \tau \in \mathbb{C} \backslash \{0\} \;.
\label{spec}
\eea
Given an $n \times n$ matrix $\M (\omega)$ with $\det \M (\omega) =1$ and $\M (\omega) = \M^{\natural} (\omega)$, we consider then
\bea
{\cal M}_{\rho,v} (\tau) = {\cal M}(\omega) \vert_{ \omega = v + \frac{\lambda}{2}    \, \rho \, \frac{\lambda -   \tau^2}{\tau} } = {\cal M} ( v + \frac{\lambda}{2}    \, \rho \, \frac{\lambda -   \tau^2}{\tau}) \;,
\label{MomMt}
\eea
which we call the monodromy matrix.

We can now state the main theorem of \cite{Aniceto:2019rhg}, where it was shown that, under very general assumptions, the canonical WH factorisation of ${\cal M}_{\rho,v} (\tau)$  w.r.t. an admissible contour $\Gamma$ in the $\tau$-complex plane, normalised at $0$, determines a solution to the field equations \eqref{eomM}.

{
\theorem
\label{theoraccr}
{\cite [Theorem 6.1]{Aniceto:2019rhg}}
Let the following assumptions hold:

(1)  
There exists an open set $S$ such that, for every 
 $(\rho_0, v_0)\in S$, one can find a
simple closed contour $\Gamma$ in the $\tau$-plane, encircling the origin and
invariant under $\tau \xmapsto[\, i_{\lambda}]{}
- \frac{\lambda}{\tau} $, such that:
\\
for all 
$(\rho,v)$ in a neighbourhood of $(\rho_0, v_0)$, 
the matrix ${\cal M}_{\rho, v} (\tau)$ given by \eqref{MomMt}, 
as well as its inverse, is 
analytic in a region $O$ 
in the $\tau$-plane containing 
$\Gamma$. 
We require $O$ to be 
 invariant under  $i_\lambda$, and such that  ${\cal M}^{\natural}_{\rho, v} (\tau)  = {\cal M}_{\rho, v} (\tau) $ on $O$;
 
 (2) for any $(\rho, v)$ in a neighbourhood of $(\rho_0, v_0)$, ${\cal M}_{\rho, v} (\tau) $ admits a canonical Wiener-Hopf factorisation
 with respect to $\Gamma$,
\bea
{\cal M}_{\rho, v} (\tau) = {\cal M}^-_{\rho, v} (\tau) \, X (\tau, \rho, v) \quad \text{on} \;\; \Gamma \;,
\label{whc}
\eea
where the ``plus" factor $X$ is normalised by
\bea
X (0, \rho, v) = \mathbb{I}_{n \times n} \;\; \text{for all $(\rho,v)$ in a region $S \subset \mathbb{R}^+ \times \mathbb{R}$} \;;
\eea

(3) the matrix function $X (\tau, \rho, v)$, for each $\tau \in \mathbb{D}^+_{\Gamma} \cup O$, and 
\bea 
M(\rho,v) : = \displaystyle{\lim_{\tau \rightarrow \infty}} {\cal M}^-_{\rho, v} (\tau) 
\label{defM}
\eea
are of class $C^2$ (w.r.t. $(\rho, v)$) and $\frac{\partial X}{\partial \rho} $ and  $\frac{\partial X}{\partial v} $
are analytic as functions of $\tau$ in $ \mathbb{D}^+_{\Gamma} \cup O$.

Then $M(\rho, v)$ defined by \eqref{defM} 
is a solution to the field equation \eqref{eomM}.\\
}

Note that these assumptions, which may at first seem complicated, allow in fact for a broad range of applicability of the results and they are easily seen to be satisfied in all known physically significant cases. Also note that different choices of the contour $\Gamma$ lead to different factorisations of the same monodromy matrix, and therefore to different solutions to the field equations originating from the same monodromy matrix \cite{Aniceto:2019rhg}.

To uniformise the notation, we will write ${\cal M}_- (\tau) $ instead of $ {\cal M}^-_{\rho, v} (\tau) $.
\\

We now apply the method discussed in Section \ref{sec:meth} to two specific $2 \times 2$ symmetric monodromy matrices $\M$ that possess a canonical WH factorisation. They can be seen as simple deformations of a monodromy matrix associated with a static solution belonging to the class of $A III$-metrics \cite{Aniceto:2019rhg}. The canonical WH factorisation of $\M$ in these two examples will however result in solutions to the vacuum field equations that, from a PDE point of view, are related in a non-trivial manner to the original solution, since deforming a solution of the vacuum field equations in such a way as to obtain another solution is a non-trivial task when employing usual PDE methods.
Deformations of this sort are sometimes called Kinnersley transformations (cf. page 4 of 
\cite{Richterek:2004bb}).

\subsection{An example with non-vanishing component $f_{2+}$}

Let 
\bea
\M (\omega) = \begin{bmatrix}
\frac{c^2}{\omega}  + s^2 \, \omega & c s \left( \frac{1}{\omega} + \omega \right)  \\
c s \left( \frac{1}{\omega} + \omega \right) & c^2 \, \omega   + \frac{s^2}{\omega} 
\end{bmatrix} = \begin{bmatrix} {\tilde a}(\omega) & {\tilde b}(\omega) \\
 {\tilde b}(\omega) & {\tilde d} (\omega) 
 \end{bmatrix} \;,
 \eea
 where $c, s \in \mathbb{C}$ with 
\bea
c^2 - s^2 = 1 \;.
\eea
Here we take $ c, s \neq 0$.

We have 
\bea
{\tilde q} (\omega) = \frac{{\tilde a}(\omega)}{{\tilde d}(\omega)} = 
\frac{ c^2 + s^2 \, \omega^2}{c^2\, \omega^2 + s^2}
\label{tqw}
\eea
and
\bea
q(\tau) = {\tilde q} \left( v + \frac{\rho}{2} \frac{1 + \tau^2}{\tau} \right) \;,
\eea
obtained from \eqref{tqw} by substituting the relation given in \eqref{spec} with $\lambda = -1$, 
\bea 
w =  v + \frac{\rho}{2} \frac{1 + \tau^2}{\tau} \;.
\eea
Note that
\bea 
w =  \frac{\rho}{2} \,  \frac{(\tau - \tau_0)}{\tau} \, ( \tau - {\tilde \tau}_0) \;,
\eea
where
\bea
\tau_0 = \frac{- v + \sqrt{v^2 - \rho^2 }}{\rho} \;\;\;,\;\;\;
{\tilde \tau}_0 =  \frac{1}{\tau_0} = \frac{- v - \sqrt{ v^2 - \rho^2 }}{\rho} \;.
\eea
We choose an admissible contour $\Gamma$ such that $\tau_0$ lies inside the contour, in which case ${\tilde \tau}_0$
lies in $\mathbb{D}^{-}_{\Gamma}$ \cite{Aniceto:2019rhg}.

The following useful relations will be used in determining a canonical factorisation of $\M_{\rho,v} (\tau)$ defined by 
\bea
\M_{\rho,v} (\tau) = \M \left( v + \frac{\rho}{2} \frac{1 + \tau^2}{\tau} \right) 
\eea
(which will simply be denoted by $\M (\tau)$ in the following):
\bea
\left( c^2 \, A^{-1} \pm s^2 \, A \right) \left( s^2 \, A^{-1} \pm c^2 \, A \right) - c^2 s^2 \left( A \pm A^{-1} \right)^2 = \pm 1
\label{pmrelA}
\eea
for any $A \in \mathbb{C}$. In particular, one obtains from \eqref{pmrelA} that
\bea
\det \M (\tau ) = 1 \;.
\eea

To establish {\it the existence of a canonical factorisation} for $\M$, we first solve the Riemann-Hilbert problem 
\bea
\M \, \phi_+ = \phi_-  \;\;\; \text{on} \;\;\; \Gamma,  \;\;\; \text{with} \;\;\; \phi_+ \in C^{\mu}_+ \;\;,\;\;  \phi_- \in C^{\mu}_{-, 0} \;\;,
\label{RHpp}
\eea
where $C^{\mu}_{-, 0}$ consists of the functions in $C^{\mu}_-$ 
vanishing at $\infty$.
Note that 
\bea 
w 
= m_- m_+ \;,
\eea
where 
\bea
m_- = \frac{\rho}{2} \, \frac{(\tau - \tau_0)}{\tau} \;\;\;,\;\;\; m_+ = \tau - {\tilde \tau}_0 \;\;\;,\;\;\; \tau_0 \in \mathbb{D}^{+}_{\Gamma} \;\;\;,\;\;\; {\tilde \tau}_0 \in \mathbb{D}^{-}_{\Gamma} \;.
\label{mmpm}
\eea
{From } \eqref{RHpp} we have that
\bea
\left\{
\begin{matrix}
\left( c^2 \, m_-^{-1} m_+^{-1} + s^2 \, m_- m_+ \right) \phi_{1+} + c s \left( m_-^{-1} m_+^{-1} + m_- m_+ \right) 
\phi_{2+} = \phi_{1-} &=& R_{1,0} \;, \\
c s \left(  m_-^{-1} m_+^{-1} + m_- m_+ \right) \phi_{1+} + \left( c^2 \, m_- m_+ + s^2 \, m_-^{-1} m_+^{-1} \right) \phi_{2+} = \phi_{2-} &=& R_{2,0} \;, 
\end{matrix}
\right.
\label{syspp12}
\eea
where $R_{1,0}$ and $R_{2,0}$ are rational functions vanishing at $\infty$ and with the same poles as $m_-, m_-^{-1}$,
i.e. with simple poles at $\tau = 0, \tau_0$. Using Cramer's rule, we have from \eqref{syspp12} that
\bea
\phi_{1+} &=& 
\begin{vmatrix}
    R_{1,0} & \quad c s \left( m_-^{-1} m_+^{-1} + m_- m_+ \right) \\
    R_{2,0} & \quad c^2 \, m_- m_+ + s^2 \, m_-^{-1} m_+^{-1}
\end{vmatrix} \nonumber\\
&=& \left( R_{1,0} \, c^2 - R_{2,0} \, c s \right) m_- m_+ + 
\left( R_{1,0}\,  s^2 - R_{2,0} \, c s \right) m_-^{-1}  m_+^{-1} \;,
\label{phi2cr}
\eea
and we look for functions $R_{1,0}, R_{2,0}$ of the form
\bea
R_{1,0} = \frac{A_1 \tau + A_0}{\tau ( \tau - \tau_0)} \;\;\;,\;\;\; 
R_{2,0} = \frac{B_1 \tau + B_0}{\tau ( \tau - \tau_0)} \;,
\eea
such that $\phi_{1+}$ given by \eqref{phi2cr} is analytic in $\mathbb{D}^{+}_{\Gamma}$, i.e. with a double zero at $0$ for the numerator of the first term in \eqref{phi2cr} and a double zero at $\tau_0$
for the second term, taking \eqref{mmpm} into account. This implies that $A_0 = B_0 = A_1 = B_1 =0$, so
\eqref{RHpp} admits only the zero solution $\phi_+ = \phi_- = 0$ and we conclude that $\M$ has a canonical factorisation.

To {\it obtain the canonical factorisation} 
\bea
\M (\tau) = \M_- (\tau) \, X(\tau)  \;\;\; \text{on} \;\;\; \Gamma,  \;\;\; \text{with} \;\;\; X(0) = \mathbb{I} \;,
\eea
we determine the two columns of $X^{-1}$ and $\M_-$ separately. The {\it first columns} of 
$X^{-1}$ and $\M_-$ are given by the (unique) solution to 
\bea
\M \, f_+ = f_- \;\;\; \text{on} \;\;\; \Gamma, \;\;\; \text{with} \;\;\; f_{\pm} \in (C^{\mu}_{\pm})^2 \;\;\;,\;\;\; f_+(0) = (1,0) \;.
\eea
Denoting $f_{\pm} = (f_{1 \pm}, f_{2 \pm} )$, we get a system as in \eqref{syspp12}, with 
$R_{1,0}$ and $R_{2,0}$ replaced by $R_1$ and $R_2$, respectively, where $R_1$ and $R_2$ are rational functions bounded at $\infty$ and with simple poles at $0$ and $\tau_0$, yielding
\bea
f_{1+} &=& 
\begin{vmatrix}
    R_{1} & \quad c s \left( m_-^{-1} m_+^{-1} + m_- m_+ \right) \\
    R_{2} & \quad c^2 \, m_- m_+ + s^2 \, m_-^{-1} m_+^{-1}
\end{vmatrix} \nonumber\\
&=& \left( R_{1} \, c^2 - R_{2} \, c s \right) m_- m_+ + 
\left( R_{1}\,  s^2 - R_{2} \, c s \right) m_-^{-1}  m_+^{-1} \;, 
\label{f1r1r2}
\\
f_{2+} &=& 
\begin{vmatrix}
c^2 \, m_-^{-1} m_+^{-1} + s^2 \, m_- m_+  & \quad R_1\\
    c s \left(  m_-^{-1} m_+^{-1} + m_- m_+ \right) & \quad R_2\\
    \end{vmatrix} \nonumber\\
&=& \left( R_{2} \, s^2 - R_{1} \, c s \right) m_- m_+ + 
\left( R_{2}\,  c^2 - R_{1} \, c s \right) m_-^{-1}  m_+^{-1} \;.
\label{phi2cr2}
\eea
Taking the form of \eqref{f1r1r2} into account, we look for $\lambda_1, \lambda_2 \in \mathbb{C}$ such that
\bea
 R_{1} \, c^2 - R_{2} \, c s  = \lambda_1 \, m_-^{-1} \;, \nonumber\\
R_{1}\,  s^2 - R_{2} \, c s = \lambda_2 \, m_- \;, 
\eea
which is equivalent to having 
\bea
R_1 &=& \lambda_1 \, m_-^{-1} - \lambda_2 \, m_- \;, \nonumber\\
R_2 &=& \frac{R_1 \, c^2 - \lambda_1 \,  m_-^{-1 } }{c s } = \frac{R_1 \, s^2 - \lambda_2 \,  m_- }{c s } 
= \frac{\lambda_1 \,  m_-^{-1} \, s^2  - \lambda_2 \,  m_- \, c^2 }{ c s } \;.
\label{r1r2l1l2}
\eea
It is easy to see that, for $R_1$ and $R_2$ given by \eqref{r1r2l1l2}, $f_{1+}$ and $f_{2+}$
given by \eqref{f1r1r2} and \eqref{phi2cr2} are analytic in $\mathbb{D}^{+}_{\Gamma} $ with
\bea
f_{1+} &=&  \lambda_1 \, m_+ + \lambda_2 \, m_+^{-1} \;, \nonumber\\
f_{2+} &=& - \lambda_2 \, \frac{c}{s} \, m_+^{-1} - \lambda_1 \, \frac{s}{c} \, m_+ \;,
\eea
and from the normalising conditions $f_{1+} (0) = 1, f_{2+} (0) = 0$, we get $\lambda_1 = c^2, \, \lambda_2 = - s^2$, so that 
\bea
f_+ = \begin{bmatrix}  c^2 \, m_+ - s^2 \, m_+^{-1} \\
s c \left( m_+^{-1} - m_+ \right) \end{bmatrix} \;\;\;,\;\;\; 
f_- = \begin{bmatrix} 
c^2 \, m_-^{-1} + s^2 \, m_- \\
s c \left( m_-^{-1} + m_- \right) \end{bmatrix} \;.
\eea

The {\it second columns in $X^{-1}$ and $\M_-$} can be determined analogously with different normalising conditions. However we will obtain them here using the results of Section \ref{sec:meth}, namely Theorem 
\ref{theo:r1r2}. Noting that, by \eqref{pmrelA},
\bea
f_{1+} \left( c^2 \, m_+^{-1} - s^2 \, m_+ \right) = 1 - f_{2 +}^2 \;,
\eea
we have that
\bea
r_1 = q \, f_{1+}^2 - f_{2+}^2 = \left( q \, f_{1+} + c^2 \, m_+^{-1} - s^2 \, m_+ \right) f_{1+} - 1 \;,
\eea
and hence 
\bea
- \left[ ( q \, f_{1+} + c^2 \, m_+^{-1} - s^2 \, m_+ ) f_{2+} \right] f_{1 +} + f_{2 +} =  r_1 \, \left( - f_{2+} \right) \;.
\label{qf1f2r1f2}
\eea
Comparing \eqref{qf1f2r1f2} with the equation for $s_{1+}$ given in \eqref{sr2pfpQQ}, we take
\bea
r_2 =  - ( q \, f_{1+} + c^2 \, m_+^{-1} - s^2 \, m_+ ) f_{2+} 
\eea
and $s_{1+} = - f_{2+}$. Using \eqref{sr2pfpQQ2} and \eqref{sr1mfmQQ2}, 
this uniquely determines $s_+$ and $s_-$ to be
\bea
s_+ = \begin{bmatrix} 
s c \left( m_+ - m_+^{-1} \right)
\\
c^2 \, m_+^{-1} - s^2 \, m_+  \end{bmatrix} \;\;\;,\;\;\; 
s_- = \begin{bmatrix} 
s c \left( m_-^{-1} + m_- \right) \\
s^2 \, m_-^{-1} + c^2 \, m_-  \end{bmatrix} \;.
\eea
In particular we see that  $s_{1+} = - f_{2+}$ and $s_{1-} = f_{2-} $.

The associated space-time metric, obtained as in \cite{Aniceto:2019rhg} from this factorisation, is given by
\bea
\Delta &=&  \frac{ \frac12 \left( v + \sqrt{v^2 -\rho^2} \right)  }{ s^2 + c^2 \frac14 \left( v + \sqrt{v^2 -\rho^2} \right)^2  } \;, \nonumber\\
B &=&  2 s\,  c  \, \rho \tau_0 = - 2 c s \left( v - \sqrt{v^2 -\rho^2} \right) \;, \nonumber\\
e^{\psi} &= & \frac{v + \sqrt{v^2 -\rho^2}}{2 \sqrt{v^2 -\rho^2} } \;.
\eea


\subsection{An example with vanishing component $f_{2+}$}

Now we consider
\bea
\M (\omega) = \begin{bmatrix}
\frac{1}{\omega}  &  \frac{\epsilon }{\omega}    \\
 \frac{\epsilon}{\omega}  &  \omega   + \frac{\epsilon^2}{\omega} 
\end{bmatrix} = \begin{bmatrix} {\tilde a}(\omega) & {\tilde b}(\omega) \\
 {\tilde b}(\omega) & {\tilde d} (\omega) 
 \end{bmatrix} \;\;\;,\;\;\; \det \M (\omega) = 1 \;,
 \eea
 where $\epsilon \in \mathbb{C}$. Here we take $\epsilon \neq 0$.

We have 
\bea
{\tilde q} (\omega) = \frac{{\tilde a}(\omega)}{{\tilde d}(\omega)} = 
\frac{ 1}{ \omega^2 + \epsilon^2}
\label{tqw2}
\eea
and
\bea
q(\tau) = {\tilde q} \left( v + \frac{\rho}{2} \frac{1 - \tau^2}{\tau} \right) \;,
\label{qeps}
\eea
obtained from \eqref{tqw2} by substituting the relation given in \eqref{spec} with $\lambda = 1$,
\bea
w =  v + \frac{\rho}{2} \frac{1 - \tau^2}{\tau} \;.
\eea
Note that
\bea 
w = - \frac{\rho}{2} \,  \frac{(\tau - \tau_0)}{\tau} \, ( \tau - {\tilde \tau}_0) \;,
\eea
where
\bea
\tau_0 = \frac{v - \sqrt{v^2 + \rho^2 }}{\rho} \;\;\;,\;\;\;
{\tilde \tau}_0 = - \frac{1}{\tau_0} = \frac{ v + \sqrt{ v^2 + \rho^2 }}{\rho} \;.
\eea
We choose an admissible contour $\Gamma$ such that $\tau_0$ lies inside the contour, in which case ${\tilde \tau}_0$
lies in $\mathbb{D}^{-}_{\Gamma}$ \cite{Aniceto:2019rhg}. 
We define
\bea
m_- =  \frac{\rho}{2} \, \tilde{\tau}_0 \, \frac{(\tau - \tau_0)}{\tau} \;\;\;,\;\;\; m_+ = - 
\frac{(\tau - {\tilde \tau}_0) }{\tilde{\tau}_0} \;\;\;,\;\;\; \tau_0 \in \mathbb{D}^{+}_{\Gamma} \;\;\;,\;\;\; {\tilde \tau}_0 \in \mathbb{D}^{-}_{\Gamma} \;.
\eea
and note that
\bea 
w 
= m_- m_+ \;.
\eea

It can be shown, as in the previous example, that the monodromy matrix 
\bea
\M_{\rho,v} (\tau) = \M \left( v + \frac{\rho}{2} \frac{1 - \tau^2}{\tau} \right) 
\eea
possesses a canonical WH factorisation \eqref{WHMf} w.r.t. $\Gamma$,  with the factors $X^{-1}$ and $\M_-$ possessing the following first columns (denoted by $f_+$ and $f_-$, respectively),
\bea
f_+ = \begin{bmatrix} f_{1+} \\ f_{2+} 
\end{bmatrix}  = \begin{bmatrix} m_+ \\ 0 
\end{bmatrix} \;\;\;,\;\;\; f_- = \begin{bmatrix} f_{1-} \\ f_{2-} 
\end{bmatrix} = \begin{bmatrix} m_-^{-1} \\  \epsilon\, m_-^{-1}
\end{bmatrix} \;.
\eea
Note that $f_{2+} (\tau) = 0$ and
$f_{2-} (\tau) = \epsilon f_{1-} (\tau)$).

To determine the second columns of $X^{-1}$ and $\M_-$, we use 
the method described in Section \ref{sec:meth}. The second columns will be denoted by $s_+$ and $s_-$, respectively,  and its components will be denoted by means of the subscripts $1, 2$.

First we express
the rational function $q$ in \eqref{qeps}  as
\bea
q = \frac{p_1 (\tau)}{p_2(\tau)} \;,
\label{qpp12}
\eea
where 
\bea
p_1 (\tau) = \tau^2 \;\;\;,\;\;\; p_2 (\tau) = \frac{\rho^2}{4} 
(\tau - \tau_0)^2 (\tau - {\tilde \tau}_0)^2 + \epsilon^2 \tau^2 \;.
\label{p1p2}
\eea
The quantity $r_1$, given in \eqref{r1qf1f2},
takes the form
\bea
r_1 = \frac{{\tilde p}_1}{p_2}
\eea
with 
\bea
{\tilde p}_1 = \frac{\tau^2 (\tau - {\tilde \tau}_0)^2}{{\tilde \tau}_0^2} 
\eea
and $p_2$ given in \eqref{p1p2}. Both ${\tilde p}_1$ and $p_2$ are polynomials of degree $4$. ${\tilde p}_1$ has two double zeroes, one located in the interior of $\Gamma$ (at $\tau =0$) and one located in the exterior of $\Gamma$ (at $\tau = {{\tilde \tau}}_0$).

Next, we determine the rational function $r_2$, by demanding that the zeroes of $(r_2 + J Q_2) f_+$ cancel the zeroes of $r_1$ in $\mathbb{D}^+_{\Gamma}$ and the zeroes of $(r_2 + J Q_1) f_-$ cancel the zeroes of $r_1$ in $\mathbb{D}^-_{\Gamma}$. 
Using \eqref{f2p0} and \eqref{ffm}, we obtain 
(where we use that $f_{2+} (\tau) = 0$ and 
$f_{2-} (\tau) = \epsilon f_{1-} (\tau)$),
\bea
\left\{ 
\begin{matrix}
R_2(0)  &=& 0 \;, \\
R'_2(0)  &=& 0 \;, \\
 R_2 ({\tilde \tau}_0)  + \epsilon p_1 ({\tilde \tau}_0)  &=& 0 \;, \\
    R'_2 ({\tilde \tau}_0) + \epsilon p'_1 ({\tilde \tau}_0)  &=& 0 \;.
 \end{matrix}
\right.
\label{4cond}
\eea
Now we use that $R_2$ is a polynomial of degree not greater than $\deg (p_2) = 4$, i.e.
\bea
R_2(\tau) = A \tau^4 + B \tau^3 + C \tau^2 + D \tau + E \;,
\eea
where the constants $A, B, C, D, E$ are determined by the four conditions \eqref{4cond} as well as by the normalising condition 
$s_{1+} (0) = 0 $ (cf. Remark \ref{normcond}).
Imposing the first two conditions in 
\eqref{4cond} gives $D = E = 0$, while imposing the two remaining conditions in \eqref{4cond} gives
\bea
B &=& - 2 A \, {\tilde \tau}_0 \;, \nonumber\\
C &=& - \epsilon + A \, {\tilde \tau}^2_0 \;.
\label{BCA}
\eea
Next we use the normalising condition 
$s_{1+} (0) = 0 $. Using the expression for $s_{1+} $ given in \eqref{sr2pfpQQ} 
we obtain $C = 0$, 
which when combined with \eqref{BCA} gives
\bea
A &=& \frac{\epsilon}{{\tilde \tau}^2_0} \;, \nonumber\\
B &=& - 2 \frac{\epsilon}{{\tilde \tau}_0 }\;.
\eea
The polynomial $R_2$ is therefore determined to be
\bea
R_2 = \epsilon \, \tau^3 \left( \frac{\tau}{{\tilde \tau}^2_0}  -  \frac{2}{{\tilde \tau}_0 } \right) \;.
\eea

Having determined $r_2$ using only knowledge about the first column vectors $f_+$ and $f_-$, the second column vectors $s_+$ and $s_-$ are uniquely determined using \eqref{sr2pfpQQ} and \eqref{sr1mfmQQ}, 
\bea
s_+ = 
\begin{bmatrix}
\epsilon \left( m_+ - m_+^{-1} \right) \\
m_+^{-1}
\end{bmatrix} \;\;\;,\;\;\;
s_-= 
\begin{bmatrix}
\epsilon \, m_-^{-1} \\
m_- + \epsilon^2 \, m_-^{-1}
\end{bmatrix} \;.
\eea

The associated space-time metric, obtained as in \cite{Aniceto:2019rhg} from this factorisation, is given by
\begin{equation}
    ds_4^2 = - \Delta \left( dt + B \, d\phi \right)^2 + \Delta^{-1} \left( e^{\psi} (  d\rho^2 + dv^2) + \rho^2 d\phi^2 \right),
\end{equation}
with 
\bea
\Delta &=&  \frac{ \frac12 \left( v + \sqrt{v^2 +\rho^2} \right)  }{ \epsilon^2 +  \frac14 
\left( v + \sqrt{v^2 +\rho^2} \right)^2  } \;, \nonumber\\
B &= & 2 \epsilon  \, \rho \tau_0 = 2 \epsilon  \left( v - \sqrt{v^2 +\rho^2} \right) \;, \nonumber\\
e^{\psi} &=& \frac{v + \sqrt{v^2 +\rho^2}}{2 \sqrt{v^2 +\rho^2} } \;.
\eea

\subsection*{Acknowledgements}
This work was partially
supported by FCT/Portugal through CAMGSD, IST-ID, projects UIDB/04459/2020 and UIDP/04459/2020.

\providecommand{\href}[2]{#2}\begingroup\raggedright

\endgroup

\end{document}